 \input amstex
\NoBlackBoxes

\documentstyle{amsppt}
\topmatter

\title{SOLVABILITY OF EQUATIONS BY QUADRATURES AND NEWTON'S THEOREM}\endtitle

\author Askold Khovanskii\endauthor

\affil{Department of Mathematics, University of Toronto, Toronto, Canada}\endaffil

\dedicatory{{dedicated to  Victor Buchstaber's 75th
birthday }  }\enddedicatory

\thanks{The work was partially supported by the Canadian Grant No. 156833-17.}\endthanks

\abstract{Picard--Vessiot theorem (1910) provides a necessary and sufficient  condition for solvability of linear differential equations of order $n$ by quadratures in terms of its Galois group. It is based on the differential Galois theory and is rather involved. J.Liouville in 1839 found an elementary criterium for such solvability for $n=2$. J.F.Ritt simplified Liouville's theorem (1948). In 1973 M. Rosenlicht proved a similar criterium  for arbitrary $n$.  Rosenlicht work relies on the valuation theory and is not elementary. In these notes we show that the elementary Liouville--Ritt method based on developing solutions in Puiseux series as functions of a parameter works smoothly for arbitrary $n$ and proves the same criterium. }\endabstract



\endtopmatter
\document
 \medskip
 \head{I. Introduction}
 \endhead
 Consider a homogeneous linear differential equation
 $$
 y^{(n)}+a_1y^{(n-1)}+\dots +a_ny=0 \tag 1
 $$
 whose coefficients  $a_i$ belong to a differential field $K$.

\proclaim {Theorem 1}
If the equation~(1)  has a non zero solution representable by generalized quadratures over $K$ then it necessarily has a solution of the  form $y_1=\exp z$ where $z'$ is  algebraic  over $K$. \footnote{ In the section 5 we  provide a generalization of Theorem 1 for  nonlinear homogeneous equations.}
\endproclaim
The following lemma is obvious.

\proclaim {Lemma 2} Assume that the equation~(1) has a non zero solution $y_1$ representable by generalized quadratures  over $K$. Then
the equation~(1)  can be solved by generalized quadratures over $K$ if and only if the linear differential equation of order $(n-1)$ over the differential field $K(y_1)$  obtained from~(1) by the reduction of order using the solution $y_1$ is solvable by generalized quadratures over  $K(y_1)$.
\endproclaim

Indeed on one hand each solution of the equation obtained from (1) by the reduction of order using $y_1$  can be expressed  in the form $(y/y_1)'$ where $y$ is a solution of~(1). On the other hand any solution $y$ of the equation $(y/y_1)'=u$, where $u$ is represented by generalized quadratures over $K(y_1)$, is representable by generalized quadratures over  $K$, assuming that $y_1$ representable by generalized quadratures  over $K$.

Thus Theorem~1 provides the following criterium for solvability of the equation~(1) by generalized quadratures.

\proclaim {Theorem 3}
The equation~(1)  is solvable by generalized quadratures over $K$ if and only if the following conditions hold:

1) the equation~(1) has a solution $y_1$ of the form $y_1=\exp  z$ where $z'=f$ is  algebraic over $K$,

2) the linear differential equation of order $(n-1)$ over  $K(y_1)$  obtained from~(1) by the reduction of order using the solution $y_1$ is solvable by generalized quadratures over  $K(y_1)$.
\endproclaim

The standard proof (E. Picard and E.~Vessiot, 1910) of Theorem 1    uses the differential Galois theory and is rather involved (see \cite 4).\footnote{A generaliztion of Theorem 1 for nonlinear homogeneous equations presented in the section 5 does not follow from the Picard-Vessio theory.}
\bigskip

In the case when the equation~(1) is a Fuchsian differential equation and $K$ is the field of rational function of one complex variable Theorem 3 has a topological explanation (see \cite 3 which allows to prove much stronger version of this result). But in general case Theorem 3 does not have a similar visual explanation.
\bigskip

In these notes I  discuss an elementary proof of Theorem~1  based on old arguments suggested by J.~Liouville, J.F.Ritt and M.~Rosenlicht.
\bigskip

Maxwell Rosenlicht in 1973 proved \cite 2  the following theorem.

\proclaim {Theorem 4}
Let $n$ be a positive integer, and let $Q$ be a polynomial in several variables with coefficients in a differential field $K$ and of total degree less than $n$. Then if the  equation
$$u^n=Q(u, u', u'', \dots) \tag 2
$$
has a solution representable by generalized quadratures over $K$, it has a solution algebraic over $K$.
\endproclaim
The logarithmic derivative $u=y'/y$ of any solution of the equation (1) satisfies the {\it associated with (1) generalized Riccati's equation of order $(n-1)$}, which is a particular case of the equation (2). Rosenlicht showed that Theorem 1 easily follows from Theorem 4 applied to the corresponding generalized Riccati's equation  (see section 5).
 In modern differential algebra abstract fields equipped with an operation of differentiation are considered. The Rosenlicht's proof of Theorem~4  is not elementary: it is applicable to abstract differential fields of characteristic zero and  makes use of the valuation theory. \footnote{ According to Michael Singer the valuation theory used in Rosenlicht's  proof  is  a fancy way of using power series methods (private communication). }.

 J. Liouvillive and J.F.Ritt deal with  fields of meromorphic functions with the usual differentiation. This point of view is natural for the analytic theory of differential equations. Solutions of differential equations are function, not elements of abstract differential fields. Analytic function could be multivalued, could have singularities and so on. For applications of results from differential algebra to the theory of differential equations some extra work is needed. In these notes we deal with functional differential fields. We present needed material in sections 3--6, mainly following the  presentation  from the book \cite 3.

 Joseph Liouville in 1839 proved Theorem~1 for $n=2$.  Joseph Fels Ritt in 1948 simplified his proof (see  [1]). The logarithmic derivative $u=y'/y$ of any solution $y$ of the homogeneous linear differential equation~(1) of second order satisfies  the  Riccati's equation
$$
u'+a_1u+a_2+u^2=0. \tag 3
$$

To prove Theorem~1 for $n=2$ J.Liouville and J.F.Ritt proved first Theorem~4 for the Riccati's equation~(3). To do that J.F.Ritt  considered a special one parametric family of solutions of (3) and  used an expansion of these solutions  as  functions of the parameter  into converging  Puiseux series. J.F.Ritt used a generalization of the following theorem based on ideas suggested by Newton.

Consider an algebraic function $z(y)$ defined by an equation  $P(y,z)=0$ where $P$ is a polynomial  with coefficients in a subfield $K$ of $\Bbb C$. Then all branches of the algebraic function $z(y)$ at the point $y=\infty$ can be developed into converging Puiseux series whose coefficients belong to a finite extension of the field $K$.

A {\it generalized Newton's Theorem} claims that the similar result holds if instead of a numerical field of coefficient one takes a field $K$ whose element are meromorphic functions on a connected Riemann surface. In the J.F.Ritt's book \cite 1 this result is proved in the same way as its classical  version using the Newton's polygon method.

Unfortunately J.F.Ritt's proof is written in  old mathematical language and  does not fit into  our presentation. Theorem 13 provides an exact statement of the generalized Newton's Theorem. It is presented  without proof: main arguments proving it are well known and classical. One also can obtain a proof modifying J.F.Ritt's exposition. Theorem 13  plays a crucial  role in these notes. For  the sake of completeness I will present its modern proof in a separate  paper.

 \bigskip

In these notes I  discuss a proof of Theorem~4 which does not rely on the valuation theory.  It generalizes J.F.Ritt's arguments (makes use of the Puiseux expansion via a generalized Newton's Theorem) and provides an elementary proof of the classical Theorem~1.

The idea of the proof goes back to Liouville and J.F.Ritt. I came up with it trying to understand  and comment the classical  book written by J.F.Ritt \cite 1.

I am grateful to Michael Singer who invited me to write comments for a new edition of this book.

\head { II. Generalized quadratures over functional differential fields}
\endhead
In the sections 1-4 we present  definitions and general statements related  to functional and abstact differential fields and classes of their extensions including extensions by generalized quadratures. We follow mainly the  presentation  from the book \cite 3. In the section 5 we define the generalized Riccati's equation and reduce Theorem 1 to Theorem 4.

\subhead{ 1. Abstract differential fields}
\endsubhead
A  field $F$  is said to be
{\it a differential field} if an additive map
$a\rightarrow a'$ is fixed that satisfies the Leibnitz rule $(ab)'=a'b+ab'$.
The element $a'$ is
called the  {\it derivative} of $a$. An element $y\in F$ is called {\it a constant} if $y'=0$.
All constants in $F$ form {\it the field of constants}. We add to the definition of  differential field an extra condition that {\it the field of constants is the field of complex numbers} (for our purpose it is enough to consider fields satisfying this condition).
An element $y\in F$  is said to be: {\it an exponential of integral} of  $a$ if $y'=ay$; {\it an integral} of $a$ if $y'=a$. In each of these cases, $y$ is defined only up to  a multiplicative  or an additive complex constant.

Let $K\subset F$ be a differential subfield in $F$. An element $y$ is said to be an {\it integral over $K$} if $y'=a\in K$. An {\it exponential of integral over $K$} is defined similarly.

Suppose that a differential field $K$ and a set $M$ lie in some differential field $F$.
{\it The adjunction} of the set $M$ to the differential field $K$ is
the minimal differential field $K\langle M\rangle$ containing both the field $K$
and the set $M$. We will refer to the transition from $K$ to $K\langle M\rangle$ as {\it adjoining} the set $M$ to the field $K$.
\definition{Definition 1}
An extension $F$ of a differential field $K$ is said to be:

1) a {\it generalized extension by integral}  if there are $y\in F$ and $f\in K$ such that $y'=f$, $y$ is transcendental  over $K$, and $F$ is a finite extension of the field $K\langle y\rangle$;

2) a {\it generalized extension by exponential of integral} if there are $y\in  F$, $f\in K$ such that $y'=fy$, $y$ is transcendental  over $K$, and  $F$ is a finite extension of the field $K\langle y\rangle$;

3) an {\it extension by generalized quadratures}
if  there exists a chain of
differential fields $K=F_0\subset \dots\subset F_n\supset F$
such that $F_{i+1}=F_i<y_i>$ for every $i=0$, $\dots$, $n-1$ where $y_i$ is an exponential of integral, an integral, or an algebraic element over the differential field $F_i$.

An element $a\in F$ is  {\it representable by generalized quadratures } over $K$, $K\subset F$, if it is contained in a certain extension of the field $K$ by generalized quadratures. The following lemma is obvious.
\enddefinition

\proclaim {Lemma 5} An extension  $K\subset F$
is an  extension by generalized quadratures
if  there is a chain  $K=F_0\subset \dots\subset F_n$
such that $F\subset F_n$ and for every $i=0$, $\dots$, $n-1$  or
$F_{i+1}$ is a finite extension of $F_i$, or
$F_{i+1}$ is  a generalized extension by integral of $F_i$, or
$F_{i+1}$ is  a generalized extension by exponential integral of $F_i$.
\endproclaim

\subhead{2. Functional differential fields and their extensions}
\endsubhead
Let $K$ be a subfield in the field $F$ of all meromorphic functions on a connected domain $U$ of the Riemann sphere $\Bbb C^1\cup \infty$ with the fixed coordinate function $x$ on $\Bbb C^1$.
Suppose that $K$ contains all complex constants and is stable under differentiation
(i.e. if $f\in K$, then $f\,'=df/dx\in K$).
Then $K$ provides an example of a {\it functional differential field.}

Let us now give a general definition.
\definition{Definition 2} Let $U,x$ be a pair consisting of a connected Riemann surface $U$ and a non constant meromorphic function $x$ on $U$. The map $f\rightarrow df/dx$ defines the derivation  in   the field $F$ of all meromorphic
functions on  $U$ (the ratio of two meromorphic 1-forms is a well-defined meromorphic function).
 A {\it functional differential field} $K$ is any differential subfield of $F$
(containing all complex constants).
\enddefinition

With any function $f$ from a functional differential field $K$ let us associate its meromorphic germ $f_a$ at  point $a\in U$. The {\it differential field $K$ is isomorphic to the differential field $K_a$ of germs at $a$ of function belonging to $K$}.

The following construction helps to {\it extend} functional differential fields.
Let $K$ be a differential subfield of the field of meromorphic functions on a
connected Riemann surface $U$
equipped with a meromorphic function $x$.
Consider any connected Riemann surface $V$ together with a
nonconstant analytic map $\pi: V\to U$.
Fix the function $\pi^*x$ on $V$.
The differential field $F$ of all meromorphic functions on $V$ with the differentiation
$\varphi'= d\varphi/\pi^*dx$ contains the differential subfield $\pi^* K$ consisting of functions
of the form $\pi^* f$, where $f\in K$.
The differential field $\pi^*K$ is isomorphic to the differential field $K$,
and it lies in the differential field $F$.
For a suitable choice of the surface $V$,
an extension of the field $\pi^*K$, which is isomorphic to $K$, can be done within
the field $F$.

Suppose that we need to extend the field $K$, say, by an integral $y$ of some function $f\in K$.
This can be done in the following way.
Consider the covering of the Riemann surface $U$ by the Riemann surface $V$ of an
indefinite integral $y$ of the form $fdx$ on the surfave $U$.
By the very definition of the Riemann surface $V$,
there exists a natural projection $\pi: V\to U$, and
the function $y$ is a single-valued meromorphic function on the surface $V$.
The differential field $F$ of meromorphic functions on $V$ with the differentiation
$\varphi'= d \varphi/\pi^*dx$ contains the element $y$ as well as the field $\pi^*K$ isomorphic to $K$.
That is why the extension $\pi^*K\langle y\rangle$ is well defined as a subfield of
the differential field $F$.
We mean this particular construction of the extension whenever we talk about extensions
of functional differential fields.
The same construction allows to adjoin
a logarithm, an exponential, an integral or an exponential of integral
of any function $f$ from a functional differential field $K$ to $K$.
Similarly, for any functions $f_1,\dots,f_n\in K$, one can adjoin a
solution $y$ of an algebraic equation $y^n +f_1y^{n-1}+\dots+f_n=0$
or all the solutions
$y_1$, $\dots$, $y_n$ of this equation to $K$
(the adjunction of all the solutions $y_1$, $\dots$, $y_n$
can be implemented on the Riemann surface of the vector-function
$\bold y= y_1$, $\dots$, $y_n$).
In the same way, for any functions $f_1,\dots,f_{n+1}\in K$, one can adjoin
the $n$-dimensional $\Bbb C$-affine space of all solutions of the
linear differential equation
$y^{(n)}+f_1y^{(n-1)}+\dots +f_ny+f_{n+1}=0$ to $K$.
(Recall that a germ of any solution of this linear differential equation admits an
analytic continuation along a path on the surface $U$ not passing through the poles
of the functions $f_1$, $\dots$, $f_{n+1}$.)

Thus, {\it all above--mentioned extensions of functional differential fields can be implemented
without leaving the class of functional differential fields}.
When talking about extensions of functional differential fields, we always mean this particular
procedure.

The differential field of all complex constants and the differential field of
all rational functions of one variable can be regarded as differential
fields of functions defined on the Riemann sphere.

\subhead {3. Classes of functions and operations on multivalued functions}
\endsubhead
An indefinite integral of an elementary function is a function rather than an element of an abstract differential field.
In functional spaces, for example, apart from differentiation and algebraic operations,
an absolutely non-algebraic operation is defined, namely, the composition.
Anyway, functional spaces provide more means for writing ``explicit formulas''
than abstract differential fields.
Besides, we should take into account that functions can be multivalued,
can have singularities and so on.

In functional spaces, it is not hard to formalize the problem of unsolvability
of equations in explicit form.
One can proceed as follows:
fix a class of functions and say that an equation is solvable explicitly if
its solution belongs to this class.
Different classes of functions correspond to different notions of solvability.

\subhead{3.1. Defining classes of functions by the lists
of data}
\endsubhead
A class of functions can be introduced by specifying a list of {\it basic functions}
and a list of {\it admissible operations}.
Given the two lists, the class of functions is defined as the set
of all functions that can be obtained
from the basic functions by repeated application of admissible operations.
Below, we define the class of functions representable by
{\it generalized quadratures} and the class of  functions representable by
{\it generalized quadratures over a functional differential field $K$} in exactly this way.

Classes of functions, which appear in the problems of solvability of differential equations by quadratures,
contain multivalued functions.
Thus the basic terminology should be made clear.
We work with multivalued functions ``globally'', which leads
to a more general understanding of classes of functions defined by
lists of basic functions and of admissible operations.
A multivalued function is regarded as a single entity.
{\it Operations on multivalued functions} can be defined.
The result of such an operation is a set of multivalued functions;
every element of this set is called a function obtained from the given functions
by the given operation.
A  {\it class of functions} is defined as the set of all (multivalued) functions
that can be obtained from the basic functions by repeated
application of admissible operations.

\subhead{3.2. Operations on multivalued functions}
\endsubhead
Let us define, for example, the sum of two multivalued functions on a connected Riemann surface $U$.

\definition{Definition 3}
Take an arbitrary point $a$ in $U$, any germ
$f_a$ of an analytic function $f$ at the point $a$ and any germ $g_a$ of an analytic function $g$
at the same point $a$.
We say that the multivalued function $\varphi$ on $U$ generated by the germ $\varphi_a=f_a+g_a$
{\it is representable as the sum of the functions} $f$ and $g$.
\enddefinition

For example, it is easy to see that exactly two functions of one variable
are representable in the form
$\sqrt{x}+\sqrt{x}$, namely, $f_1=2\sqrt{x}$ and $f_2\equiv 0$.
Other operations on multivalued functions are defined in exactly the same way.
{\it For a class of multivalued functions, being stable under addition means that,
together with any pair
of its functions, this class contains all functions representable as their sum.}
The same applies to all other operations on multivalued functions understood
in the same sense as above.

In the definition given above, not only the operation of addition plays a key role but
also the operation of analytic continuation hidden in the notion
of multivalued function.
Indeed, consider the following example.
Let $f_1$ be an analytic function defined on an open subset $V$ of the complex line
$\Bbb C^1$ and admitting no analytic continuation outside of $V$,
and let $f_2$ be an analytic function on
$V$ given by the formula $f_2=-f_1$.
According to our definition, the zero function is representable in the form
$f_1+f_2$ {\it on the entire complex line}.
By the commonly accepted viewpoint, the equality $f_1+f_2=0$ holds inside
the region $V$ but not outside.

Working with multivalued functions globally, we do not insist on the existence of
{\it a common region}, were all necessary operations would be performed on
single-valued branches of multivalued functions.
A first operation can be performed in a first region,
then a second operation can be performed in a second, different region
on analytic continuations of functions obtained on the first step.
In essence, this more general understanding of operations is equivalent to including
analytic continuation to the list of admissible operations on the analytic germs.
\subhead{4. Functions representable by generalized quadratures}
\endsubhead
In this section we define  functions of one complex variable representable by generalized quadratures and functions representable by generalized quadratures  over a functional differential field. We also discuss a relation of these notions with extensions of functional differential fields by generalized quadratures. First we'll present needed lists of basic functions and of admissible operations.

\proclaim {List of basic elementary functions}

1. All complex constants and an independent variable $x$.

2. The exponential, the logarithm, and the power $x^\alpha$ where $\alpha$ is any constant.

3. The trigonometric functions sine, cosine, tangent, cotangent.

4. The inverse trigonometric functions arcsine, arccosine, arctangent, arccotangent.
\endproclaim
\proclaim{Lemma 6}
Basic elementary functions can be expressed through the exponentials and
the logarithms with the help of complex constants, arithmetic operations
and compositions.
\endproclaim
Lemma 6 can be considered as a  simple exercise. Its proof can be found in [3].

\proclaim {List of some classical  operations}

1)  Operation of composition takes functions $f$,$g$ to the function $f\circ g$.

2)  The arithmetic operations take functions $f$, $g$ to the functions $f+g$, $f-g$, $fg$, and $f/g$.

3) The operation of differentiation takes function $f$ to the function $f'$.

4) The operation of integration takes function $f$ to a solution of equation $y'=f$ (the function $y$ is defined up to an additive constant).

5) The operation of taking  exponential of integral takes function $f$ to a solution of equation $y'=fy$ (the function $y$ is defined up to a multiplicative constant).

6)  The operation of solving algebraic equations takes functions $f_1,\dots,f_n$ to the function $y$ such that $y^n+f_1y^{n-1}+\dots+f_n=0$ ( the function $y$ is not quite uniquely determined by functions $f_1,\dots,f_n$ since an algebraic equation of degree $n$ can have $n$ solutions).

\endproclaim

\definition {Definition 4} The class of functions of one complex variable {\it representable by generalized quadratures} is defined by the following data:

List of basic functions: basic elementary functions.

List of admissible operations: Compositions, Arithmetic operations,  Differentiation, Integration,Operation of taking exponential of integral, Operation of solving algebraic equations.
\enddefinition
\proclaim{Theorem 7}
A (possibly multivalued) function of one complex variable  belongs to the class of functions representable by generalized quadratures  if and only if it belongs to some extension   of the differential field of all constant functions of one variable by generalized quadratures.
\endproclaim

Theorem 7 follows from Lemma 6 (All needed  arguments can be found in [3]).

Let $K$ be  a functional differential field  consisting of meromorphic functional on a connected Riemann surface $U$ equipped with a meromorphic function $x$.

\definition {Definition 5} The class of functions representable by {\it generalized quadratures over the functional differential field $K$} is defined by the following data:

List of basic functions: all functions from the field $K$.

List of admissible operations: Operation of composition with a function $\phi$ representable by generalized quadratures that takes $f$ to $\phi\circ f$, Arithmetic operations,  Differentiation, Integration, Operation of taking exponential of integral, Operation of solving algebraic equations.
\enddefinition

\proclaim{Theorem 8}
A (possibly multivalued) function on the Riemann surface $U$ belongs to the class of generalized quadratures over a functional differential field $K$ if and only if it belongs to some  extension of  $K$ by generalized quadratures.
\endproclaim

Theorem 8 follows from Lemma 6 (all needed  arguments can be found in [3]).

\subhead{5. Generalized Riccati's equation}
\endsubhead
In this section we define the generalized Riccati's equation and reduce Theorem 1 to Theorem 4. In this section we also generalize Theorem~1 for nonlinear  homogeneous equations (this generalization will not be used in the next sections).

Assume that $u$ is the logarithmical derivative of a non identically equal to zero meromorphic function $y$, i.e the relation $y'=uy$ holds.
\definition {Definition 6} Let $D_n$ be a polynomial in $u$ and in its derivatives $u,u',\dots,u^{(n-1)}$ up to order $(n-1)$ defined by induction by
the following conditions:
$$
D_0=1; \ \  D_{k+1}=\frac{d D_k}{d x}+u D_k.
$$
\enddefinition
\proclaim {Lemma 9} 1) The polynomial $D_n$ has integral coefficients and $\deg D_n= n$. The degree $n$ homogeneous part of $D_n$ equals to $u^n$ (i.e. $D_n=u^n+\tilde D_n$ where $\deg \tilde D_n<n$). 2) If $y$ is a function whose logarithmic derivative equals to $u$ (i.e. if $y'=uy$) then for any $n\geq 0$ the relation $y^{(n)}=D_n(u) y$ holds.
\endproclaim
Both claims of Lemma 9 can be easily checked by induction.

Consider a homogeneous linear differential equation (1)
whose coefficients  $a_i$ belong to a differential field $K$.

 \definition{Definition 7} The  equation
  $$
 D_n+a_1D_{n-1}+\dots +a_nD_0=0 \tag 4
 $$
of  order $n-1$ is called the {\it generalized Riccati's equation} for the homogeneous linear differential equation (1).

 \enddefinition

 \proclaim {Lemma 10} A non identically equal to zero function $y$ satisfies the linear differential equation (1) if and only if its logarithmic derivative $u=y'/y$ satisfies the generalized Riccati's equation (4).
 \endproclaim
 \demo{Proof} Let $y$ be a nonzero solution of (1) and let $u$ be its logarithmic derivative. Then dividing (1) by $y$ and using the identity $y^{(k)}/y=D_k(u)$ we obtain that $u$ satisfies (4). If $u$ is a solution of (4) then multiplying (4) by $y$ and using the identity $y^{(k)}=D_k(u)y$ we obtain that $y$ is a non zero solution of (1).

 \enddemo
\proclaim {Corollary 11} 1) The equation (1) has a non zero solution  representable by generalized quadratures  over $K$ if and only if the equation (4) has a solution  representable by generalized quadratures  over $K$.

2) The equation (1) has a solution $y$ of the form $y=\exp z$ where $z'=f$ is an algebraic function over $K$ if and only if the equation (4) has an algebraic solution over $K$.
\endproclaim

\demo{Proof} 1) A non zero function $y$ is representable by generalized quadratures over $K$ if and only if  its logarithmic derivative $u=y'/y$ is representable by generalized quadratures over $K$.

2) A function $y$ is equal to $\exp z$ where $z'=f$ if and only if its logarithmic derivative is equal to $f$.
\enddemo

The generalized Riccati's equation (4) satisfies the conditions of Theorem 4. Thus Theorem 1 follows from Theorem 4 and from  corollary 11.

Let us generalize the results of this section. Consider an order $n$ homogeneous equation
 $$
P(y,y',\dots,  y^{(n)})=0 \tag 1'
 $$
 where $P$ is a degree $m$ homogeneous polynomial in $n+1$ variables $x_0,x_1,\dots, x_{n}$ over a functional differential field $K$.
\definition{Definition 7'} The  equation
  $$
P(D_0,D_1, \dots, D_n)=0 \tag 4'
 $$
of  order $n-1$  is called the {\it generalized Riccati's equation} for the  homogeneous equation (1').
\enddefinition

\proclaim {Lemma 10'} A non identically equal to zero function $y$ satisfies the homogeneous  equation (1') if and only if its logarithmic derivative $u=y'/y$ satisfies the generalized Riccati's equation (4').
 \endproclaim

\proclaim {Corollary 11'} 1) The equation (1') has a non zero solution  representable by generalized quadratures  over $K$ if and only if the equation (4') has a solution  representable by generalized quadratures  over $K$.

2) The equation (1') has a solution $y$ of the form $y=\exp z$ where $z'=f$ is an algebraic function over $K$ if and only if the equation (4') has an algebraic solution over $K$.
\endproclaim

Lemma 10' and Corollary 11' can be proved exactly in the way as Lemma 10 and Corollary 11'.

Let us defined the {\it $\xi$-weighted degree $\deg_\xi x^p$} of the monomial  $x^p=x_0^{p_0}\cdot\dots\cdot x_n^{p_n}$ by the following formula: $$ \deg_\xi x^p=
\sum_{i=0}^{i=n} i m_i
.$$
We will say that a polynomial $P(x_0,\dots,x_n)$ satisfies the {\it $\xi$-weighted degree condition} if the sum of coefficients of all monomials in $P$ having the biggest $\xi$-weighted degree is not equal to zero. A polynomial $P$ having a unique  monomial with the biggest $\xi$-weighted degree  automatically satisfies this condition. For example a degree $m$ polynomial $P$ containing a term $ax_n^m$ with $a\neq 0$  automatically satisfies $\xi$-weighted degree condition.

\proclaim {Theorem 1'} Consider the homogeneous equation~(1') with the polynomial $P$ satisfying the $\xi$-weighted degree condition. If this equation  has a non zero solution representable by generalized quadratures over $K$ then it necessarily has a solution of the  form $y_1=\exp z$ where $z'$ is  algebraic  over $K$.
\endproclaim
\demo{Proof} It is easy to check that if the polynomial $P$ satisfies the $\xi$-weighted degree condition then the generalized Riccati's equation (4') satisfies the conditions of Theorem 4. Thus Theorem 1' follows from  Theorem 4 and corollary 11'.
\enddemo

\remark{Remark} There exists a complete analog of Galois theory for linear homogeneous differential equations (see [4]). Theorem 1 can be proved using this theory. The differential Galois group of a nonlinear homogeneous differential equation (1') could be very small and for such equation a complete analog of Galois theory  does not exist. Thus Theorem 1' can not be proved in a similar way.
\endremark

\head {III. Special extensions of functional differential fields}
\endhead

In the chapter we consider simple extensions of functional differential fields. We also consider algebraic extensions  of fields of rational functions over functional fields. In the section 6  we state  Theorem 13 (generalized Newton's Theorem) playing a crucial role for these notes. Section 9.3 and 10.3 contains first steps for  our proof of  Theorem 4.

\subhead {6. Finite extensions of fields of rational functions}
\endsubhead
In this section we will discuss finite extensions of the field $K(y)$ of rational functions over a subfield $K$ of the field of meromorphic function on a connected Riemann surface $U$.

Let $F$ be  extension of $K(y)$ by a root $z$ of a degree $m$ polynomial $P(z)\in (K[y])[z]$ over the ring $K[y]$ irreducible over the field $K(y)$. Let $X$ be the product $U\times \Bbb C^1$ where $\Bbb C^1$ is the standard complex line with the coordinate function $y$. An element of the field $K(y)$ can be considered as meromorphic function on $X$. One can associate with the element $z\in F$  a multivalued  algebroid  function on $X$ defined by equation $P(z)=0$. Let
$D(y)$ be the discriminant of the polynomial $P$. Let $\Sigma\subset U\times \Bbb C^1=X$ be the hypersurface defined by equation $p_m (y)\cdot D(y)=0$ where $p_m(y)$ is the leading coefficient of the polynomial $P$.

\proclaim{Lemma 12} 1) About a point $x\in X\setminus \Sigma$ the equation $P(z)=0$ defines $m$  germs $z_i$ of analytic functions whose values at  $x$ are simple roots of polynomial $P$. 2) Let $x$ be the point $(a,y)\in U\times \Bbb C^1\setminus \Sigma$. Then the field $F$ is isomorphic to the extension $K_a(y,z_i)$ of the field $K_a$ of  germs at $a\in U$ of  functions from the field $K$ (considered as germs at $x=(a,y)$ of functions independent of $y$) extended by the independent variable $y$ and by the germ $z_i$ at $x$ satisfying the equation $P(z)=0$.
\endproclaim
\demo {Proof} The statement 1) follows from the implicit function theorem. The statement 2) follows from 1).
\enddemo

Below we state Theorem 13 which is a generalization of Newton's Theorem  about expansion of an algebraic functions into converging Puiseux series. It is stated  without proof  (see comments in the introduction).

We use notations introduced in the beginning of this section. Let $z_i$ be  an element  satisfying a polynomial equation $P(z)=0$ over the ring $K[y]$,  where $K$ is a subfield  of the field of meromorphic functions on $U$. Then there exists {\it a finite extension $K_P$ of the field $K$ associated with the polynomial $P$} such that the following theorem holds.

\proclaim {Theorem 13} There is a finite covering $\pi :U_P\rightarrow U\setminus O_P$ where $O_P\subset U$ is a discrete subset, such that  the following properties hold:

1) the extension $K_P$ can be realized by a subfield of the field of meromorphic functions on $U_P$ containing the field $\pi^*K$ isomorphic to $K$.

2) there is a continuous positive function $r:U\setminus O_P\rightarrow \Bbb R$ such that in the  open domain $W\subset (U\setminus O_P)\times \Bbb C^1$ defined by the inequality $|y|>r(a)$  all $m$ germs  $z_i$ of $z$ at a point $(a,y)$ can be developed into converging   Puiseux series

$$
z_i=z_{i_k}  y^{\frac {k}{p}}+z_{i_{k-1}}y^{\frac {k-1}{p}}+\dots \tag 5
$$
whose coefficients $z_{i_j}$ are germs of analytic functions at the point $a\in U\setminus O_P$ having analytical  continuation as regular functions on  $U_P$  belonging to the field  $K_P$.
\endproclaim

Theorem 13 for a special  case when the field $K$ is a subfield of the field of complex numbers it is natural to name the Newton's Theorem.  One can consider
$K$ as a field of constant functions on any connected Riemann surface $U$.  One can chose $O_P$ to be the empty set, $U_P$ to be equal $U$, projection  $\pi:U_P\rightarrow U$ to be  the identity map, the function $r:U\rightarrow \Bbb R$  to be a big enough constant. In this case Theorem 13 states that an algebraic function $z$ has a Puiseux expansion at infinity whose coefficients belong to a finite extension $K_P$ of the field $K$. This statement can be proved by Newton's polygon method.

Let $F$ be an extension of $K(y)$ by a root $z$ of polynomial $P$ and let $K_P$ be the finite extension of the field $K$ introduced in Theorem 13. The extension $F_P$ of the field $K_P(y)$ by $z$ is easy to deal with. Denote the product $U_P\times \Bbb C^1$ by $X_P$.

\proclaim {Lemma 14} Let $x\in X_P$ be the point $(a,y_0)\in U_P\times \Bbb C^1$. Then the field $F_P$ is isomorphic to the extension $K_{P,a}(y,z_i)$ of the field $K_{P,a}$ of  germs at $a\in U_P$ of  functions from the field $K_P$ (considered as germs at $x=(a,y_0)$ of functions independent of $y$) extended by the independent variable $y$ and by the germ  at $x$  of the function $z_i$ defined by (5).
\endproclaim

Lemma 14 follows from Theorem 13.

\subhead  {7.  Finite extensions of differential fields}
\endsubhead
In this section  we  discuss finite extensions of functional differential fields.

Let
$$
P(z)=z^n+a_1z^{n-1}+\dots+a_n
$$
be an
irreducible polynomial over  $K$,
$P\in K[z]$.
Suppose that a functional differential field $F$ contains $K$ and a root $z$ of $P$.

\proclaim{Lemma 15}
 The field $K(z)$ is stable under the differentiation.

\endproclaim
\demo{Proof} Since  $P$ is irreducible over $K$, the polynomial
  $\frac{\partial P}{\partial z}$ has no common roots with $P$ and is different from zero
  in the field $K[z]/(P)$. Let $M$ be a polynomial  satisfying a congruence $M\frac{\partial P}{\partial z}\equiv -
  \frac{\partial P}{\partial x}\pmod P$.
  Differentiating the identity $P(z)=0$ in the field $F$,
we obtain that $\frac{\partial P}{\partial z}(z)z'+
\frac{\partial P}{\partial x}(z)=0$, which implies that $z'=M(z)$.
Thus the derivative of the element $z$ coincides with the value at $z$ of
the polynomial $M$. Lemma 15 follows from this fact.

\enddemo

Let $K\subset F$ and $\hat K\subset \hat F$ be  functional differential fields, and
$P$, $\hat P$
irreducible polynomials over $K$,
$\hat K$ correspondingly.
Suppose that  $F$, $\hat F$ contain  roots $z$, $\hat z$
of  $P$, $\hat P$.

\proclaim {Theorem 16} Assume that there is an isomorphism  $  \tau:K\rightarrow \hat K$ of differential fields $K$, $\hat K$ which maps coefficients of the polynomial $P$ to the corresponding coefficients of the polynomial $\hat P$. Then $\tau$ can be extended in a unique way to the differential isomorphism $ \rho:K(z)\rightarrow \hat K(\hat z)$.
\endproclaim

Proof of  Theorem 16 could be obtain by  the  arguments  used in the proof of Lemma 15.

\subhead {8.  Extension by one  transcendental element}
\endsubhead
Let $U$ be a connected Riemann surface and let $K$ be a differential field of meromorphic functions on $U$. Let $\Bbb C^1$ be the standard complex  line with the coordinate function $y$. Elements of the field $K(y)$ of rational functions over $K$ could be considered as meromorphic  functions on $X=U\times \Bbb C^1$.

In the field $K(y)$   there are two natural operations of differentiations. The first operation $R(y)\rightarrow \frac{\partial R}{\partial x}(y)$ is defined as follows: the derivative $\frac{\partial }{\partial x}$ of the independent variable $y$ is equal to zero, and derivative $\frac{\partial }{\partial x}$ of an element $a\in K$ is  equal to its derivative $a'$ in the field $K$. For the second operation $R(y)\rightarrow \frac{\partial R}{\partial y}(y)$ the derivative of an element $a\in K$  is equal to zero and the derivative of the independent variable $y$ is equal to one.

Let $K\subset F$  be  differential fields and let $\theta\in F$  be a transcendental  element over $K$. Assume that  $\theta '\in K\langle \theta\rangle$. Under this assumption the field $K\langle\theta\rangle$ has a following description.

\proclaim{Lemma 17} 1) The map  $\tau:K\langle \theta \rangle \rightarrow K(y)$ such that $\tau (\theta)=y$ and $\tau(a)=a$ for $a\in K$ provides an isomorphism between the field  $K\langle \theta \rangle$ considered  without the operation of differentiation and the field $K(y)$ of rational functions over $K$. 2) If $\tau (\theta')=w\in K(y)$ then for any $R\in K(y)$ and $z\in K\langle \theta \rangle $ such that $\tau(z)=R$ the following identity holds
$$
\tau (z')=\frac{\partial R}{\partial x}+ \frac{\partial R}{\partial y}w. \tag 7
$$
\endproclaim
\demo{Proof}The first claim of the lemma is straightforward. The second claim follows from the chain rule.
\enddemo

Let $\Theta\subset X=U\times \Bbb C^1$  be the graph of function $\theta:U\rightarrow \Bbb C^1$. The following lemma is straightforward.

\proclaim{Lemma 18} The differential field $K\langle\theta \rangle$ is isomorphic to the field  $K(y)|_\Theta$ obtained by restriction on $\Theta$ of functions from the field $K(y)$ equipped with the differentiation given by (7). For any point $a\in \Theta$ The differential field $K\langle\theta \rangle$ is isomorphic to the differential field of germs at $a\in \Theta$ of functions from $K(y)|_\Theta$.
\endproclaim

\subhead{9. An  extension by integral}
\endsubhead
In this section we consider  extensions of transcendental degree one of a differential field $K$ containing an integral $y$ over $K$ which does not  belong to $K$, $y\notin K$.
\subhead {9.1. A pure transcendental  extension by integral}
\endsubhead
Let $\theta$ be an integral over $K$,i.e $\theta'=f\in K$. Assume that $\theta$ is a transcendental element over $K$. \footnote{Easy to check that if an integral $\theta$ over $K$ does not belong to $K$, then $\theta$ is a transcendental element over $K$ (see [3]). We will not use this fact.}

\proclaim {Lemma 19} 1) The field $K\langle \theta \rangle$ is isomorphic to the field $K(y)$ of rational functions over $K$ equipped with the following differentiation
$$
R'=\frac{\partial R}{\partial x}+ \frac{\partial R}{\partial y}f.\tag 8
$$
2) For every complex number $\rho\in \Bbb C$ the map $\theta\rightarrow \theta+\rho$ can be extended to the unique isomorphism $G_\rho:K\langle \theta \rangle\rightarrow K\langle \theta \rangle$ which fixes elements of the field $K$.

3) Each isomorphism of $K\langle \theta \rangle$ over $K$ is an isomorphism $G_\rho$ for some $\rho\in \Bbb C$. Thus the Galois group of $K\langle \theta \rangle$ over $K$ is the additive group of complex numbers $\Bbb C$.
\endproclaim

\demo{Proof} The claim 1) follows from Lemma 17. For any $\rho \in \Bbb C$ the element $\theta_\rho=\theta+\rho$ is a transcendental element over $K$ and $\theta_\rho'$ equals to $f$.  Thus the claims 2) is correct.  The claim 3) followst from 2) because if $y'=f$ then $y=\theta_\rho$ for some $\rho\in \Bbb C$.
\enddemo

\subhead {9.2.  A generalized extension by integral}
\endsubhead
According to Lemma 19 the differential field $K\langle\theta\rangle$ is isomorphic to the field $K(y)$ with the differentiation given by (8). Let $F$ be an extension of $K\langle\theta\rangle$ by an  element $z\in F$   which satisfies some equation $\tilde P(z)=0$ where $\tilde P$ is an irreducible polynomial over $K\langle \theta \rangle$. The isomorphism between $K\langle\theta\rangle$ and $K(y)$ transforms the polynomial $\tilde P$ into some polynomial $P$ over $K(y)$. Below we  use notation from section 6 and deal we the multivalued algebroid function $z$ on $X$  defined by $P(z)=0$.

Assume that at a point $x\in X$ there are germs of analytic functions $z_i$ satisfying the equation $P(z_i)=0$. Let $\theta_\rho$ be the function $(\theta+\rho):U_P\rightarrow \Bbb C^1$ and let $\Theta_\rho \subset X=U\times\Bbb C^1$ be its graph. The point  $x=(p,q)\in U\times \Bbb C^1$ belongs to the graph $\Theta_{\rho(x)}$ for $\rho(x)=q-\theta(p)$.

Let $K(y)|_{\Theta_{\rho(x)}}$ be the differential field of germs at the point $x\in \Theta_{\rho(x)}$ of  restrictions on $\Theta_{\rho(x)}$ of functions from the field $K(y)$ equipped with the differentiation given by (8).

\proclaim{Lemma 20} The differential field $F$ is isomorphic to the finite extension of the differential field  $K(y)|_{\Theta_{\rho(x)}}$ obtained by adjoining the germ at $x\in \Theta_{\rho(x)}$ of the restriction to $\Theta_{\rho(x)}$ of an analytic  germ $z_i$ satisfying~$P(z_i)=0$.
\endproclaim
\demo{Proof} For the trivial extension  $F=K\langle \theta \rangle$  Lemma 20 follows from  Lemmas 18, 19. Theorem 16 allows to complete the proof for non trivial finite extensions $F$ of $K\langle \theta \rangle$.
\enddemo

According to section 6 with the polynomial $P$ over $K(y)$ one can associate the finite extension $K_P$ of the field $K$ and the Riemann surface $U_P$ such that Theorem 13 holds. Since $K$ is  functional differential field the field $K_P$ has a natural structure of  functional differential field. Below we will apply Lemma 20 taking instead of $K$ the field $K_P$ and considering the extension $F_P\supset K_P\langle \theta \rangle$ by the same algebraic element $z\in F$. The use of $K_P$ instead of $K$ allows to apply the expansion (5) for~$z_i$.

\proclaim {Theorem 21} Let $x\in X_P= U_P\times \Bbb C^1$ be a point $(a,y_0)$ with $|y|>>0$. The differential field $F_P$ is isomorphic to the extension of the differential field  of  germs at the point $a\in U_P$ of functions from the differential fiels $K_P$ by the  following germs: by the germ at $a$ of the integral $\theta_{\rho(x)}$ of the function $f\in K$, and by a germ at $a$ of the composition $z_i(\theta_\rho)$ where $z_i$ is a germ at $x$ of a function given by a Puiseux series (5).
\endproclaim

\demo{Proof} Theorem 21 follow from Lemma 20 and Theorem 13.
\enddemo

\subhead {9.3. Solutions of equations in a generalized extension by integral}
\endsubhead
We will use notations from sec. 9.1 and 9.2.

Let $T(u,u',\dots,u^{(n)}$ by a polynomial in independent function $u$ and its derivatives with coefficients from the functional differential field $K$. Consider the equation
$$
T(u,u',\dots,u^{(N)})=0. \tag 9
$$
In general  the derivative of the highest order $u^{(N)}$ cannot  be represented as a function  of other derivatives via the relation (9). Thus even existence of local solutions of (9) is problematic and we have no information about global behavior of its solutions.

Assume that (9) has a solution $z$ in a generalized extension by integral $F\supset K\langle \theta \rangle$ of $K$. The solution $z$ has a nice global property: it is a meromorphic function on a Riemann surface $U_P$ with a projection $\pi:U_P\rightarrow U$ which proves a locally trivial covering above $U\setminus O_P$, where $O_P\subset U$ is discrete subset.

Moreover existence of a solution $z$ implies existence of a  family $z(\rho)$ of similar solutions  depending on a parameter $\rho$: one  obtains such family of solutions  by using an integral $\theta+\rho$ instead of the integral $\theta$ (see Lemma 20). If the parameter $\rho$ has big absolute value $|\rho|>>0$ for a point $a\in U_P$ of the germ $z(\rho)$ can be expand in the Piueux series in $\theta_\rho$:

$$
z_i(\rho)=z_{i_k}  \theta_\rho^{\frac {k}{p}}+z_{i_{k-1}}\theta_\rho^{\frac {k-1}{p}}+\dots \tag 10
$$

The series is converging and so it can be differentiate using the relation $\theta_\rho'=f$.

\proclaim{Lemma 22} If $z_{i_k}'\neq 0$ then the leading term of the Puiseux series for $z_i(\rho)'$ is $z_{i_k}'  \theta_\rho^{\frac {k}{p}}$. Otherwise the leading term has degree $<\frac {k}{p}$. The leading term of the derivative of any order of $z_{i_k}$ has degree $\leq \frac {k}{p}$.
\endproclaim

Let us plug into the differential polynomial $T(u,u',\dots,u^{(N)}) $ the germ (10) and develop the result into Puiseux series in $\theta_\rho$.  If the germ  $z_i(\rho)$ is a solution of the equation (9) then all terms of this Puiseux series are equal to zero. In particular the leading coefficient is zero. This observation is an important step for  proofing Theorem~4.

\subhead{10. An  extension by an exponential of integral}
\endsubhead
In this section we consider  extensions of transcendental degree one of a differential field $K$ containing an exponential integral $y$ over $K$ which is not algebraic over  $K$.
\subhead {10.1. A pure transcendental  extension by an exponential integral}
\endsubhead
Let $\theta$ be an an exponential of integral over $K$, i.e $\theta'=f\theta$ where $f \in K$. Assume that $\theta$ is a transcendental element over $K$. \footnote{Easy to check that if an an exponential of integral $\theta$ over $K$ is algebraic over  $K$, then $\theta$ is a radical over $K$, i.e. $\theta^k\in K$ for some positive integral $k$ (see [3]). We will not use this fact.}

\proclaim {Lemma 23} 1) The field $K\langle \theta \rangle$ is isomorphic to the field $K(y)$ of rational functions over $K$ equipped with the following differentiation
$$
R'=\frac{\partial R}{\partial x}+ \frac{\partial R}{\partial y}fy. \tag 11
$$
2) For every complex number $\mu\in \Bbb C^*$ not equal to zero  the map $\theta\rightarrow \mu \theta$ can be extended to the unique isomorphism $G_\mu:K\langle \theta \rangle\rightarrow K\langle \theta \rangle$ which fixes elements of the field $K$.

3) Each isomorphism of $K\langle \theta \rangle$ over $K$ is an isomorphism $G_\mu$ for some $\mu\in \Bbb C^*$. Thus the Galois group of $K\langle \theta \rangle$ over $K$ is the multiplicative group of complex numbers $\Bbb C^*$.
\endproclaim

\demo{Proof} The claim 1) follows from Lemma 17. For any $\mu \in \Bbb C^*$ the element $\theta_\mu=\mu \theta$ is a transcendental element over $K$ and $\theta_\mu'=f\theta$.  Thus the claims 2) is correct.  The claim 3) follows from 2) because if $y'=fy$ and $y\neq 0$ then $y=\theta_\mu$ for some $\mu\in \Bbb C^*$.
\enddemo

\subhead {10.2.  A generalized extension by exponential of integral}
\endsubhead
According to Lemma 23 the differential field $K\langle\theta\rangle$ is isomorphic to the field $K(y)$ with the differentiation given by (11). Let $F$ be an extension of $K\langle\theta\rangle$ by an  element $z\in F$  which satisfies some equation $\tilde P(z)=0$ where $\tilde P$ is an irreducible polynomial over $K\langle \theta \rangle$. The isomorphism between $K\langle\theta\rangle$ and $K(y)$ transforms the polynomial $\tilde P$ into some polynomial $P$ over $K(y)$. Below we  use notation from section 6 and deal we the multivalued algebroid function $z$ on $X$  defined by $P(z)=0$.

Assume that at a point $x\in X$ there are germs of analytic function $z_i$ satisfying the equation $P(z_i)=0$. Let $\theta_\mu$ be the function $(\mu \theta):U_P\rightarrow \Bbb C^1$ and let $\Theta_\mu \subset X=U\times\Bbb C^1$ be its graph. The point  $x=(p,q)\in U\times \Bbb C^1$ where $q\neq 0$ belongs to the graph $\Theta_{\mu(x)}$ for $\mu(x)=q\cdot\theta(p)^{-1}$.

Let $K(y)|_{\Theta_{\mu(x)}}$ be the differential field of germs at the point $x\in \Theta_{\mu(x)}$ of  restrictions on $\Theta_{\mu(x)}$ of functions from the field $K(y)$ equipped with the differentiation given by (11).

\proclaim{Lemma 24} The differential field $F$ is isomorphic to the finite extension of the differential field  $K(y)|_{\Theta_{\mu(x)}}$ obtained by adjoining the germ at $x\in \Theta_{\mu(x)}$ of the restriction to $\Theta_{\mu(x)}$ of an analytic  germ $z_i$ satisfying $P(z_i)=0$.
\endproclaim
\demo{Proof} For the trivial extension  $F=K\langle \theta \rangle$  Lemma 24 follows from  Lemmas 18, 23. Theorem 16 allows to complete the proof for non trivial finite extensions $F$ of $K\langle \theta \rangle$.
\enddemo

According to section 6 with the polynomial $P$ over $K(y)$ one can associate the finite extension $K_P$ of the field $K$ and the Riemann surface $U_P$ such that Theorem 13 holds. Since $K$ is  functional differential field the field $K_P$ has a natural structure of  functional differential field. Below we will apply Lemma 24 taking instead of $K$ the field $K_P$ and considering the extension $F_P\supset K_P\langle \theta \rangle$ by the same algebraic element $z\in F$. The use of $K_P$ instead of $K$ allows to apply the expansion (5) for~$z_i$.

\proclaim {Theorem 25} Let $x\in X_P= U_P\times \Bbb C^1$ be a point $(a,y_0)$ with $|y|>>0$. The differential field $F_P$ is isomorphic to the extension of the differential field  of  germs at the point $a\in U_P$exponential of integral $\theta_{\mu(x)}$, where $\theta_{\mu(x)}'=f \theta_{\mu(x)}$ for the function $f\in K$, and by a germ at $a$ of the composition $z_i(\theta_\mu)$ where $z_i$ is a germ at $x$ of a function given by a Puiseux series (5).
\endproclaim
\demo{Proof} Theorem 25 follow from Lemma 24 and Theorem 13.
\enddemo

\subhead {10.3. Solutions of equations in a generalized extension by exponential of integral}
\endsubhead
Assume that (9) has a solution $z$ in a generalized extension by exponential of integral $F\supset K\langle \theta \rangle$ of $K$. Solution $z$ has a nice global property: it is a meromorphic function on a Riemann surface $U_P$ with a projection $\pi:U_P\rightarrow U$ which proves a locally trivial covering above $U\setminus O_P$, where $O_P\subset U$ is discrete subset.

Moreover existence of a solution $z$ implies existence of  family $z(\mu)$ of similar solutions  depending on a parameter $\mu\in \Bbb C^*$: one  obtains such family of solutions  by using an exponential of integral $\mu \theta$ instead of the exponential of integral $\theta$ (see Lemma 24). If the parameter $\mu$ has big absolute value $\mu|>>0$ for a point $a\in U_P$ of the germ $z(\mu)$ can be expand in the Piueux series in $\theta_\mu$:

$$
z_i(\mu)=z_{i_k}  \theta_\mu^{\frac {k}{p}}+z_{i_{k-1}}\theta_\mu^{\frac {k-1}{p}}+\dots \tag 12
$$

The series is converging and so it can be differentiate using the relation $\theta_\mu'=f\theta_\mu$.

\proclaim{Lemma 26} If $z_{i_k}'+\frac{k}{p}z_{i_k}\neq 0$ then the leading term of the Piueux series for $z_i(\mu)'$ is $(z_{i_k}'+\frac{k}{p}z_{i_k})  \theta_\mu^{\frac {k}{p}}$. Otherwise the leading term has degree $<\frac {k}{p}$. The leading term of derivative of any order of $z_{i_k}$ has degree $\leq \frac {k}{p}$.
\endproclaim

Let us plug into the differential polynomial $T(u,u',\dots,u^{(N)}) $ the germ  (12) and develop the result into Puiseux series in $\theta_\mu$.  If the germ  $z_i(\mu)$ is a solution of the equation (9) then all terms of this Piueux series are equal to zero. In particular the leading coefficient is zero. This observation is an important step for  proofing Theorem~4.

\head {IV. Proof of Rosenlicht's Theorem}
\endhead

Here we  complete an elementary  proof of Theorem 4 discovered by Maxwell Rozenlicht \cite 2.
We will proof first the simplier theorems 27, 28 of a similar   nature.

\proclaim{Theorem 27}
Assume that the equation (2) over a functional differential field $K$ has a solution $z\in F$ where $F$ is a generalized extension by integral  of  $K$. Then (2) has a solution in the algebraic extension $K_P$  of $K$ associated with the polynomial $P$.
\endproclaim

\demo{Proof} If the constant term of the differential polynomial $T(u,u',u'',\dots)= u^n-Q(u, u', u'', \dots)$ is equal to zero, then (2) has solution $u\equiv 0$ belonging to $K$. In this case we have nothing to prove.

Below we will assume that the constant term $T_0$ of $T$ is not equal to zero.
Thus the differential polynomial $T$ has two special terms: the term $u^n$ which is the only term of highest degree $n$ and the term $T_0$ which is the only term of smallest degree zero.

Assume that (2) has a solution $z$ in a generalized extension by integral $F\supset K\langle \theta \rangle$ of $K$. According to sec. 9.3 the existence of such a solution $z$  implies the existence of  family $z(\rho)$ of germs of solutions  depending on a parameter $\rho$ such that when the whose   absolute value $\rho|$ is big enough $z(\rho)$ can be expanded in  Puiseux series (10)  in $\theta_\rho$.

We will show that the degree $\frac{k}{p}$ of the leading term in (10) is equal to zero and the leading coefficient $z_{i_0}\in K_P$ satisfies (2). This will proof Theorem 27.

According to Lemma 22 the  leading term of the derivative of any order of $z_{i_k}$ has degree $\leq \frac {k}{p}$. Thus the leading term of Puiseux series obtained by plugging (10) instead of $u$ into differential polynomial $Q$ has degree $< n \frac {k}{p}$. The leading term of the Puiseux series obtained by arising (10 ) to the  $n$-th power is equal to $n \frac {k}{p}$. If $\frac {k}{p}>0$ this term can not be canceled after plugging (10) instead of $u$ into differential polynomial $T$. Thus the degree $\frac {k}{p}$ can not be positive.

Let us plug (10) into the differential polynomial $T(u,u',\dots)-T_0 $. We will obtain a Puiseux series of negative degree if $\frac {k}{p}<0$. Thus the term $T_0$ in the sum $(T-T_0)I+T_0$ can not be canceled. Thus $\frac {k}{p}$ can not be negative.

We proved that $\frac {k}{p}=0$. If in this case we plug (10) into the differential polynomial $T(u,u',\dots)$ we obtain a Puiseux series having only one term of nonnegative degree which is equal to zero.
From Lemma 22 it is easy to see that this term equals to $T(z_{i_0},z_{i_0}',\dots) \theta_\rho^0$. Thus $z_{i_0}\in K_P$ is a solution of (2). Theorem 27 is proved.
\enddemo
\proclaim{Theorem 28}
Assume that  equation (2) over a functional differential field $K$ has a solution $z\in F$ where $F$ is a generalized extension by an exponential of integral  of  $K$. Then (2) has a solution in the algebraic extension $K_P$  of $K$ associated with the element $z\in F$.
\endproclaim
\demo {Proof} Theorem 28 can be proved exactly in the same way as Theorem 27. Just instead of Lemma 22 one has to use  Lemma 26. In the case when  leading term of the Puiseux expansion of $z_\mu$ has degree zero, the leading coefficient of its derivative equals to $z_{i_0}'$ (see Lemma 26). That is why the case $\frac {k}{p}=0$ in Theorem 28 can be treated exactly in the same way as in Theorem 27.
\enddemo
Now we ready to prove  Theorem 4.

\demo {Proof of Theorem 4} By assumption the equation (2)  has a solution $z\in F$ where $F$ is an extension of $K$ by generalized quadratures. By Lemma 5 there is a chain  $K=F_0\subset \dots\subset F_m$ such that $F\subset F_m$ and for every $i=0$, $\dots$, $m-1$  or $F_{i+1}$ is a finite extension of $F_i$, or $F_{i+1}$ is  a generalized extension by integral of $F_i$, or $F_{i+1}$ is  a generalized extension by exponential integral of $F_i$. We prove Theorem 4 by induction in the length $m$ of the chain of extension. For $m=1$ Theorem 4 follows from Theorem 27, or from Theorem 28. Assume that  Theorem 4 is true for $m=k$. A chain $F_0 \subset F_1\subset \dots\subset F_{k+1}$ provides the chain $F_1\subset \dots\subset F_{k+1}$ of extensions of length $k$ for the field $F_1$. Thus (2) has an algebraic solution $z$ over the field $F_1$.
The extension $F_0\subset \tilde F_1$, where $\tilde F_1$ is the extension of $F_1$ by the element $z$, is either an algebraic extension, or extension by generalized integral or extension by generalized exponential of integral. Thus for the extension $F_0\subset \tilde F_1$ Theorem 4 holds. We completed the inductive proof of Theorem 4.

\enddemo

\end